\documentclass[12pt]{amsproc}

\usepackage[]{graphics,graphicx,float,latexsym,epsfig}
\usepackage{amsfonts,amstext,amsmath,amssymb,times}

\newcommand{\ignore}[1]{} 

\usepackage[lmargin=1in,rmargin=1in,tmargin=1in,bmargin=1in]{geometry}
\numberwithin{equation}{section} 

\newtheorem{theorem}{Theorem}
\newtheorem{lemma}{Lemma}
\newtheorem{corollary}{Corollary}

\theoremstyle{definition} 

\newtheorem{remark}{Remark}

\newcommand{\de}{\delta}
\newcommand{\De}{\Delta}
\newcommand{\La}{\Lambda}
\newcommand{\la}{\lambda}
\newcommand{\si}{\sigma}

\newcommand{\al}{\alpha}
\newcommand{\te}{\theta}
\newcommand{\ga}{\gamma}

\newcommand{\vae}{\varepsilon}
\newcommand{\tiX}{\widetilde{X}}
\newcommand{\tite}{\tilde{\te}}
\newcommand{\vano}{\varnothing}
\newcommand{\alo}{\al \to 0}

\newcommand{\mrm}[1]{\mathrm{#1}}

\newcommand{\ADD}{{\mrm{ADD}}}

\newcommand{\D}{{\mrm{D}}}
\newcommand{\Exp}{{\mrm{Exp}}}

\newcommand{\mb}[1]{\mathbf{#1}} 
\newcommand{\Pb}{{\mb{P}}}
\newcommand{\Eb}{{\mb{E}}}
\newcommand{\Xb}{{\mb{X}}}

\newcommand{\Fb}{{\mb{F}}}
\newcommand{\Ibo}{{\mb{I}}}
\newcommand{\Sib}{{\mb{\Sigma}}}

\newcommand{\mbs}[1]{\boldsymbol{#1}} 
\newcommand{\Dal}{{\mbs{\Delta}}_\infty(\al)}
\newcommand{\Deb}{{\mbs{\Delta}}}

\newcommand{\nub}{{\mbs{\nu}}}

\def\One{\mathchoice{\rm 1\mskip-4.2mu l}{\rm 1\mskip-4.2mu l}
{\rm 1\mskip-4.6mu l}{\rm 1\mskip-5.2mu l}}
\newcommand\Ind[1]{{\One_{\{#1\}}}}

\newcommand{\mc}[1]{\mathcal{#1}} 

\newcommand{\Nc}{{\mc{N}}}
\newcommand{\Fc}{{\mc{F}}}

\newcommand{\xra}{\xrightarrow} 
\newcommand{\xla}{\xleftarrow}

\newcommand{\bl}{\bigl}
\newcommand{\br}{\bigr}
\newcommand{\Bl}{\Bigl}
\newcommand{\Br}{\Bigr}

\newcommand{\set}[1]{\left\{#1\right\}}

\newcommand{\brc}[1]{\left(#1\right)}

\renewcommand{\le}{\leqslant} 
\renewcommand{\ge}{\geqslant}

\font\msbmx=msbm10                   
\font\msbmvii=msbm7                  
\font\msbmv=msbm5
\def\varnothing{\mathchoice{\mbox{\msbmx\char'077}}%
{\mbox{\msbmx\char'077}}{\mbox{\msbmvii\char'077}}{\mbox{\msbmv\char'077}}}%

\def\esup{\mathop{\rm ess\,sup}\limits}

\begin{document}

\vspace*{-0.3in}

\centerline{\small \textsl{Submitted to the Theory of Probability and Its Applications, August
2006}}

\bigskip

\title[Asymptotic optimality in Bayesian change-point detection problems ]
{\normalsize \bf ASYMPTOTIC OPTIMALITY IN BAYESIAN CHANGE-POINT DETECTION PROBLEMS UNDER GLOBAL
FALSE ALARM PROBABILITY CONSTRAINT}

\author[A.G.\ Tartakovsky] {{\sc Alexander G.\ Tartakovsky}\\
{\em University of Southern California\\
 Department of Mathematics\\
Los Angeles, CA 90089-2532, USA}}

\address{
\parbox{3.5in} {Department of Mathematics \\
University of Southern California \\
3620 S.\ Vermont Ave, KAP-108\\
Los Angeles, CA 90089-2532, USA} }

\email{tartakov@usc.edu}


\date{\today}


\begin{abstract}

In 1960s Shiryaev developed Bayesian theory of change detection in independent and identically
distributed (i.i.d.) sequences. In Shiryaev's classical setting the goal is to minimize an average
detection delay under the constraint imposed on the average probability of false alarm. Recently,
Tartakovsky and Veeravalli (2005) developed a general Bayesian asymptotic change-point detection
theory (in the classical setting) that is not limited to a restrictive i.i.d.\ assumption. It was
proved that Shiryaev's detection procedure is asymptotically optimal under traditional average
false alarm probability constraint, assuming that this probability is small. In the present paper,
we consider a less conventional approach where the constraint is imposed on the global, supremum
false alarm probability. An asymptotically optimal Bayesian change detection procedure is proposed
and thoroughly evaluated for both i.i.d.\ and non-i.i.d. models when the global false alarm
probability approaches zero.

\

\emph{Keywords and Phrases:} Bayesian change-point detection, sequential detection, asymptotic
optimality, global false alarm probability, nonlinear renewal theory, non-i.i.d.\ observations,
$r$-quick convergence.

\end{abstract}

\maketitle

%

%
\section{Introduction} \label{sec:intro}

The classical change-point detection problem deals with the i.i.d.\ case where there
is a sequence of observations $X_1,X_2,\dots$ that  are identically distributed with
a probability density function (pdf) $f_0(x)$ for $n < \la$ and with a pdf $f_1(x)$
for $n \ge \la$, where $\la$, $\la =1,2,\dots$ is an unknown point of change. In
other words, the joint pdf of the vector $\Xb_1^n=(X_1,\dots,X_n)$ conditioned on
$\la=k$ has the form
\begin{equation} \label{eq:pdfiid}
p(\Xb_1^n|\la=k)= \begin{cases} \prod_{i=1}^{k-1} f_0(X_i) \times \prod_{i=k}^{n}
f_1(X_i), & \text{if $k \le n$} \\
\prod_{i=1}^{n} f_0(X_i), & \text{if $k > n$}.
\end{cases}
\end{equation}

More generally, the observations may be nonidentically distributed or correlated or both, i.e.,
non-i.i.d. In the most general non-i.i.d.\ case the model can be described as follows
\begin{equation} \label{eq:pdfnoniid}
p(\Xb_1^n|\la=k)= \begin{cases} \prod_{i=1}^{k-1} f_{0}(X_i|\Xb_1^{i-1}) \times
\prod_{i=k}^{n} f_{1}(X_i |\Xb_1^{i-1}), & \text{if $k \le n$} \\
\prod_{i=1}^{n} f_{0}(X_i|\Xb_1^{i-1}), & \text{if $k > n$},
\end{cases}
\end{equation}
where $f_{0}(X_i|\Xb_1^{i-1})$ and $f_{1}(X_i|\Xb_1^{i-1})$ are conditional densities for $X_i$
given $\Xb_1^{i-1}=(X_1,\dots,X_{i-1})$ that may depend on $i$. In addition, the post-change pdf
$f_1(X_i|\Xb_1^{i-1})$ may depend on the point of change $k$.

A change-point detection procedure $\tau$ is a stopping time with respect to the sequence of
sigma-algebras $\Fc_n=\si(\Xb_1^n)$, $n \ge 1$, i.e., $\{\tau \le n\} \in \Fc_n$, $n\ge 0$.

Let, for any $\la=k < \infty$, $\Pb_{k}$ ($\Eb_k$) be the probability measure (expectation) under
which the conditional pdf of $X_n$ is $f_{0}(X_n|\Xb^{n-1})$ if $n \le k-1$ and is
$f_{1}(X_n|\Xb^{n-1})$ if $n \ge k$. If $\la=\infty$, i.e., when the change does not occur,
$\Pb_\infty$ ($\Eb_\infty$) is the probability measure (expectation) under which the conditional
pdf of $X_n$ given $\Xb^{n-1}$ is $f_{0}(X_n|\Xb^{n-1})$ for every $n \ge 1$.

For $\la=k$, a true detection happens when $\tau \ge k$ and false if $\tau < k$. The design of the
quickest change detection procedures involves optimizing the tradeoff between a ``risk" $R_k(\tau)$
related to the detection delay $(\tau-k)^+$ and a loss $L_k(\tau)$ due to a false alarm. Possible
risk functions are $R_k(\tau)=\Eb_k(\tau-k|\tau \ge k)$ and
$R_k(\tau)=\esup\Eb_k[(\tau-k)^+|\Fc_{k-1}]$. The first one was introduced by
Pollak~\cite{PollakAS85} and the second one by Lorden~\cite{LordenCPDAS71}. The loss $L_k(\tau)$
can be measured by the mean time to false alarm $\Eb_k [\tau \Ind{\tau < k}]$ or by the probability
of false alarm (PFA) $\Pb_k(\tau <k)$. Note that since $\set{\tau <k}\in \Fc_{k-1}$,
\begin{equation}
\Pb_k(\tau <k)=\Pb_\infty(\tau <k) \quad \text{and} \quad \Eb_k[ \tau \Ind{\tau < k}]=\Eb_\infty
[\tau \Ind{\tau < k}].
\end{equation}
Therefore, the requirements of controlling the PFA $\Pb_k(\tau <k)$ and the mean time to false
alarm $\Eb_k [\tau \Ind{\tau < k}]$ for all $k \ge 1$ are equivalent to controlling
$\sup_k\Pb_\infty(\tau <k)=\Pb_\infty(\tau < \infty)$ and $\Eb_\infty \tau$, respectively. Note
that the requirement of having $\Pb_\infty(\tau < \infty) \le \al$, $\al <1$ leads to $\Eb_\infty
\tau=\infty$ and the requirement $\Eb_\infty \tau=\ga$, $\ga < \infty$ leads to $\Pb_\infty(\tau <
\infty) =1$.

Under the constraint on the mean time to false alarm $\Eb_\infty \tau \ge \ga$, $\ga
>0$, a uniformly optimal detection procedure that minimizes the average detection
delay $\Eb_k(\tau-k|\tau \ge k)$ or $\esup\Eb_k[(\tau-k)^+|\Fc_{k-1}]$ for all $k \ge 1$ does not
exist and one has to resort to the minimax setting of minimizing $\sup_k R_k(\tau)$. In the i.i.d.\
case, Lorden \cite{LordenCPDAS71} showed that the CUSUM detection test is asymptotically optimal
with respect to the essential supremum speed of detection measure
$\sup_k\esup\Eb_k[(\tau-k)^+|\Fc_{k-1}]$ for low false alarm rate as $\ga\to\infty$. Later,
Moustakides \cite{MoustakidesAS86} improved this result showing that the CUSUM test is actually
exactly optimal for all $\ga >0$ if the threshold can be chosen in such a way that $\Eb_\infty \tau
= \ga$. See also Ritov \cite{Ritov90} for an alternative proof of this property. More recently,
Shiryaev \cite{ShiryaevRMS96} and Beibel \cite{Beibel96} proved the same result for the problem of
detecting a change in the mean value of a continuous-time Brownian motion. Pollak \cite{PollakAS85}
introduced the randomized at the initial point Shiryaev-Roberts test, which will be referred to as
the Shiryaev-Roberts-Pollak (SRP) test, and proved that this test is nearly optimal with respect to
$\sup_k\Eb_k(\tau-k|\tau \ge k)$ as $\ga \to \infty$. Further, Lai \cite{LaiIEEE98} and Tartakovsky
\cite{TarPreprint03} proved that these both detection tests are asymptotically (first order)
optimal as $\ga \to \infty$ for fairly general non-i.i.d.\ models. More recently, Fuh
\cite{Fuh03,Fuh04} proved asymptotic optimality of the CUSUM and SRP procedures for hidden Markov
models.

Specifically, let $Z_n^k$ denote the log-likelihood ratio between the hypotheses ``$H_k:\la=k$" and
$H_\infty$,
\begin{equation} \label{LLRgen}
Z_n^k = \sum_{i=k}^n \log \frac{f_{1}(X_i|\Xb^{i-1})}{f_{0}(X_i|\Xb^{i-1})}, \quad k \le n,
\end{equation}
and assume that $(n-k)^{-1} Z_n^k \to q$ almost surely (a.s.) as $n \to\infty$ under $\Pb_k$, where
$q$ is a positive and finite number. Assuming in addition a certain rate of convergence in the
above strong law, it follows from \cite{LaiIEEE98,TarPreprint03} that
\[
\liminf_{\ga\to\infty} \frac{\inf_{\{\tau: \Eb_\infty \tau \ge \ga\}}
\sup_k\Eb_k(\tau-k|\tau \ge k)}{\log \ga} \ge 1/q,
\]
which is attained for CUSUM and SRP tests with the threshold $h=\log \ga$.

Further generalizations to composite hypotheses, nonparametric problems, multipopulation problems,
multisensor distributed change detection problems, as well as detailed discussions of several
challenging application areas were presented in Tartakovsky
\cite{TarSteklov94,TartakovskyIEEECDC05}, Tartakovsky et al \cite{TartakovskyetalStamet06}, and
Tartakovsky and Veeravalli \cite{TarVeerASMbook04}.

On the other hand, for the standard CUSUM and SRP tests (with constant thresholds), the ``global"
PFA $\Pb_\infty(\tau<\infty)=1$. To guarantee the condition $\Pb_\infty(\tau<\infty) \le \al$ for
$\al<1$ in these latter tests, one may use a curved stopping boundary that increases in time in
place of the constant threshold. Borovkov~\cite{BorovkovTVP98} proved that the CUSUM and SRP tests
with certain curved thresholds are asymptotically optimal for i.i.d.\ data models with respect to
the conditional average detection delay (ADD) $\Eb_k(\tau-k | \tau \ge k)$ as $k\to\infty$. It
follows from the latter work that when $k$ is large, the conditional ADD of these procedures
increases as $O(\log k)$. This happens because of the very strong supremum probability constraint.
Therefore, under this constraint neither minimax nor uniform solutions are feasible in asymptotic
setting when $\al\to 0$, since for any small $\al$ there exists a large $k$ that cannot be
neglected. We argue that under the constraint imposed on the global (supremum) PFA the only
feasible solution is Bayesian. Indeed, in the Bayesian setting, due to averaging the increasing
threshold generates a constant term that can be neglected when $\al$ is small.

If, however, the false alarm rate is measured in terms of the local PFA $\sup_k \Pb_\infty(k \le
\tau \le k +T-1)$ or by the local conditional PFA $\sup_k \Pb_\infty(k \le \tau \le k +T-1| \tau
\ge k)$ in some time-window $T$, which may go to infinity at a certain rate, then the CUSUM and SRP
detection tests have uniformly asymptotically optimal properties, i.e., minimize the conditional
ADD $\Eb_k(\tau-k | \tau \ge k)$ for every $k \ge 1$ (cf.\ Lai \cite{LaiJRSS95,LaiIEEE98} and
Tartakovsky \cite{TartakovskyIEEECDC05}).

In Shiryaev's classical Bayesian setting (see Shiryaev \cite{ShiryaevSMD61}-\cite{Shiryaevbook78}
and Peskir and Shiryaev \cite{PeskirShiryaev-book06}), there is a prior distribution
$\pi_k=\Pb(\la=k)$, $k\ge 0$, and the constraint is imposed on the average false alarm probability
\[
\Pb^\pi(\tau<\la) = \sum_{k=1}^\infty \pi_k\Pb_k(\tau < k),
\]
i.e., $\Pb^\pi(\tau<\la) \le \al$, $\al <1$. The goal is to find an optimal procedure that
minimizes the average detection delay
\[
\Eb^\pi (\tau - \la)^+ = \sum_{k=0}^\infty \pi_k \Eb_k(\tau-k)^+
\]
in the totality of procedures $\{\tau: \Pb^\pi(\tau<\la) \le \al\}$ or an asymptotically optimal
procedure that minimizes the delay when $\al\to 0$ (see Tartakovsky and Veeravalli
\cite{TartakovskyVeeravalli-TVP2004} and Baron and Tartakovsky \cite{BaronTartakovskySA05}). Here
$\Pb^\pi$ ($\Eb^\pi$) is the average probability measure (expectation) defined as
$\Pb^\pi(\Omega)=\sum_{k=0}^\infty \pi_k \Pb_k(\Omega)$.

Shiryaev \cite{Shiryaevbook78} proved that the stopping time
\begin{equation} \label{Shiryaevst}
\nu_B= \min\set{n: \Pb(\la \le n | \Fc_n) \ge B}
\end{equation}
is optimal in the i.i.d.\ case and for the geometric prior distribution if the threshold is chosen
so that $\Pb^\pi(\nu_B<\la)=\al$. Yakir \cite{Yakir94} generalized this result for Markov models.
Recently, Tartakovsky and Veeravalli \cite{TartakovskyVeeravalli-TVP2004} and Baron and Tartakovsky
\cite{BaronTartakovskySA05} proved that the Shiryaev stopping time with the threshold $B_\al=1-\al$
is asymptotically optimal as $\al\to 0$ for a wide class of prior distributions and non-i.i.d.\
models under very general conditions. Moreover, it follows from
\cite{BaronTartakovskySA05,TartakovskyVeeravalli-TVP2004} that the Shiryaev detection test
minimizes (asymptotically) not only the average detection delay $\Eb^\pi (\tau - \la)^+$ but also
higher positive moments of the detection delay $\Eb^\pi[(\tau-\la)^m|\tau \ge \la]$, $m \ge 1$.

Note once again that the event $\set{\tau<k}$ belongs to the sigma-field
$\Fc_{k-1}=\si(\Xb_1^{k-1})$, which implies $\Pb_k(\tau<k)=\Pb_\infty(\tau<k)$.
Therefore,
\begin{equation} \label{PktoPinfty}
\Pb^\pi(\tau < \la)=\sum_{k=1}^\infty\pi_k \Pb_\infty(\tau<k).
\end{equation}

Another possibility is to impose a more strong, supremum constraint
\[
\sup_{k \ge 1} \Pb_k(\tau <k)= \sup_{k \ge 1} \Pb_\infty(\tau < k)= \Pb_\infty (\tau < \infty) \le
\al,
\]
i.e., to consider the class of stopping times $\Dal=\{\tau: \Pb_\infty (\tau < \infty) \le \al\}$
for which the worst-case (global) false alarm probability $\sup_{k \ge 1} \Pb_k(\tau <k)$ is
restricted by the given number $\al <1$. The goal is to find an optimal procedure from the
following optimization problem
\[
\inf_{\tau \in \Dal} \Eb^\pi (\tau-\la)^+ \rightarrow \tau_{\mrm{opt}}.
\]
As we already mentioned above, the minimax solution is not feasible under this strong constraint~--
the minimax delay is infinitely large. We believe that the only feasible solution is Bayesian.
However, see Assaf et al \cite{AssafetalAS93} and Remark 1 in Section \ref{s:Remarks} regarding a
dynamic sampling technique in minimax problems.

In this paper, we are interested in the latter optimization problem. However, it is difficult to
find an exact solution to this optimization problem even in the i.i.d.\ case. For this reason, we
focus on the asymptotic problem, letting $\al$ go to zero. Since $\Eb^\pi (\tau-\la)^+/\Pb^\pi(\tau
\ge \la)$ and, by \eqref{PktoPinfty}, $\Pb^\pi(\tau \ge \la)\ge 1-\Pb_\infty (\tau < \infty)$, this
latter asymptotic problem is equivalent to minimizing the average detection delay (ADD) of the form
\[
\inf_{\tau \in \Dal} \Eb^\pi (\tau-\la | \tau \ge \la) \quad \text{as $\al \to 0$}.
\]
Moreover, we will address the problem of minimizing higher moments of the detection
delay
\[
\inf_{\tau \in \Dal} \Eb^\pi [(\tau-\la)^m | \tau \ge \la], \; m > 1 \quad \text{as
$\al \to 0$}.
\]
We will write $\ADD^\pi(\tau)= \Eb^\pi (\tau-\la | \tau \ge \la)$ and $\D_m^\pi(\tau)= \Eb^\pi
[(\tau-\la)^m | \tau \ge \la]$ for brevity.

Beibel~\cite{Beibel00} considered a purely Bayesian problem for the Brownian motion
with the risk function $c \Eb^\pi (\tau-\la)^+ + \Pb_\infty (\tau < \infty)$ when
the cost of detection delay $c$ goes to zero and the loss due to the false alarm is
measured by $\Pb_\infty (\tau < \infty)$.

In the present paper, we show that the techniques developed in \cite{BaronTartakovskySA05,
DragTarVeerIEEEIT99,LaiIEEE98,TarSISP98,TartakovskyVeeravalli-TVP2004} can be effectively used for
studying asymptotic properties of change-point detection tests in the class $\Deb_\infty(\al)$ when
$\al\to 0$ for general stochastic models.

\section{The Detection Procedure} \label{s:DetProc}

Let $``H_k: \la =k"$ and $``H_\infty: \la=\infty"$ denote the hypotheses that the change occurs at
the point $\la=k$ ($k < \infty$) and does not occur. The likelihood ratio between these hypotheses
based on the observation vector $\Xb^n = (X_1,\dots,X_n)$ is
\[
\La_n^k := \frac{p(\Xb^n|\la=k)}{p(\Xb^n|\la=\infty)} = \prod_{i=k}^n
\frac{f_{1}(X_{i}|\Xb^{i-1})}{f_{0}(X_{i}|\Xb^{i-1})}, \quad n \ge k
\]
(see \eqref{eq:pdfnoniid}).

We will always use the convention that for $n=0$ , i.e., before the observations become available,
$\La_0^0= f_{1}(X_0)/f_{0}(X_0)=1$  almost everywhere. For the sake of convenience and with very
little loss of generality, we will also assume that $\pi_0=0$. Since $\La_0^0=1$, the likelihood
ratios $\La_n^0$ and $\La_n^1$ are equal, which means that the hypotheses $\la=0$ and $\la=1$ are
not distinguishable and, therefore, introducing a positive mass at the point $\la=0$ has little
practical meaning. Generalization to the case where $\pi_0>0$ is straightforward.

Define the statistic
\[
G_n = \sum_{k=1}^\infty \pi_k \prod_{i=k}^n \frac{f_{1}(X_{i}|\Xb^{i-1})}{f_{0}(X_{i}|\Xb^{i-1})},
\quad G_0=1,
\]
which is nothing but the average likelihood ratio of the hypotheses $H_k$ and $H_\infty$, and
introduce the stopping time
\begin{equation} \label{eq:tauA}
\tau_A=\min\set{n\ge 1: G_n \ge A}, \quad A > 1.
\end{equation}
Note that the statistic $G_n$ can be represented in the following form
\begin{equation} \label{eq:Gn}
G_n = \sum_{k=1}^n \pi_k e^{Z_n^k} + \Pi_{n+1}, \quad n \ge 0,
\end{equation}
where $\Pi_{n+1}=\Pb (\la \ge n+1)$ and $Z_n^k = \log \La_n^k$ is the log-likelihood ratio (LLR)
between the hypotheses $H_k$ and $H_\infty$ given in \eqref{LLRgen}.

It is useful to establish a relationship between the detection procedure $\tau_A$ and Shiryaev's
stopping time $\nu_B$ defined in \eqref{Shiryaevst}. Making use of the Bayes formula and
\eqref{eq:Gn}, we obtain
\begin{align*}
\Pb(\la \le n | \Xb_1^n) & = \frac{\sum_{k=1}^n \pi_k p(\Xb_1^n|\la=k)}{\sum_{k=1}^\infty \pi_k
p(\Xb_1^n|\la=k)} = \frac{\sum_{k=1}^n \pi_k \La_n^k}{\sum_{k=1}^\infty \pi_k \La_n^k} = \frac{G_n
- \Pi_{n+1}}{G_n},
\end{align*}
which shows that the stopping time $\tau_A$ can be written as
\[
\tau_A = \min \set{n\ge 1: \Pb(\la \le n |\Xb_1^n)\ge 1- \Pi_{n+1}/A}, \quad A >1.
\]

Therefore, while in Shiryaev's test the posterior probability $\Pb(\la \le n|\Xb_1^n)$ is compared
to a constant threshold, in the proposed detection test the threshold is an increasing function in
$n$. This is an unavoidable penalty for the very strong supremum PFA constraint.

\section{The Upper Bound on the Global Probability of False Alarm} \label{s:PFA}

Let $\Pb(\Fc_n)$ denote the restriction of the measure $\Pb$ to the $\si$-algebra
$\Fc_n=\si(X_1,\dots,X_n)$. The following lemma gives a simple upper bound for the PFA
$\Pb_\infty(\tau_A<\infty)$ in a general case. This conservative bound will be improved in Section
\ref{sss:HO}, Lemma \ref{le3} in the i.i.d.\ case.

\begin{lemma} \label{le1} 
For any $A>1$,
\begin{equation} \label{boundPFA}
\Pb_\infty(\tau_A < \infty) \le 1/A.
\end{equation}
\end{lemma}


\proof

Noting that
\[
G_n = \frac{d\Pb^\pi(\Fc_n)}{d\Pb_\infty(\Fc_n)}
\]
and using the Wald likelihood ratio identity, we obtain
\[
\Pb_\infty(\tau_A < \infty) = \Eb_\infty \Ind{\tau_A<\infty} = \Eb^\pi[
G_{\tau_A}^{-1}\Ind{\tau_A<\infty}] .
\]
By definition of the stopping time $\tau_A$, the value of $G_{\tau_A} \ge A$ on the
set $\{\tau_A<\infty\}$, which implies inequality \eqref{boundPFA}.

\endproof


Therefore, setting $A=A_\al=1/\al$ guarantees $\Pb_\infty(\tau_A < \infty) \le \al$,
i.e.,
\[
A_\al=1/\al \Rightarrow \tau_{A_\al} \in \Dal.
\]

\section{Asymptotic Optimality and Asymptotic Performance} \label{s:AO}

\subsection{The asymptotic lower bound for moments of the detection delay} \label{ss:lowerbound}

The proof of asymptotic optimality of the detection procedure $\tau_A$ with $A=A_\al=1/\al$ as $\al
\to 0$ is performed in two steps. The first step is to obtain an asymptotic lower bound for moments
of the detection delay $\D_m^\pi(\tau)$ for any procedure from the class $\Dal$. The second step is
to show that the procedure $\tau_{A_\al}$ achieves this lower bound.

It turns out that the second step is case dependent. For example, proofs and corresponding
conditions of asymptotic optimality are different in the i.i.d.\ and non-i.i.d.\ cases. See
Remark~\ref{rem1} in Section \ref{ss:AOnoniid}. For this reason, we will consider these two cases
separately. However, for deriving the lower bound the same techniques can be used in all cases. We
start with deriving the lower bound in a general, non-i.i.d.\ case.

Define $L_\al=q^{-1}|\log \al|$, $L_A= q^{-1}\log A$  and, for $0<\vae <1$,
\begin{align*}
\ga_{\vae,\al}^{\pi}(\tau) & =\Pb^\pi\set{\la \le \tau < \la+(1-\vae)L_\al},
\\
\ga_{\vae,A}^{\pi}(\tau_A) & =\Pb^\pi\set{\la \le \tau_A < \la+(1-\vae)L_A},
\end{align*}
where $q$ is a positive finite number.

The number $q$ plays a key role in the asymptotic theory. In the general case, we do not specify
any particular model for the observations. As a result, the LLR process has no specific structure.
We hence have to impose some conditions on the behavior of the LLR process at least for a large
$n$. It is natural to assume that there exists a positive finite number $q=q(f_1,f_0)$ such that
$n^{-1}Z_{k+n-1}^k$ converges almost surely to $q$, i.e.,
\begin{equation}\label{asgeneral}
\frac{1}{n}Z_{k+n-1}^{k} \xra[n\to\infty]{\Pb_k-\text{a.s.}} q \quad \text{for every $k<\infty$}.
\end{equation}
As we discuss in the end of this section, \eqref{asgeneral} holds in the i.i.d.\ case with
$q=I=\Eb_1 Z_1^1$ whenever the Kullback-Leibler information number $I$ is positive and finite.
Therefore, in the general case the number $q$ plays the role of the Kullback-Leibler number, and it
can be treated as the asymptotic local divergence of the pre-change and post-change models
(hypotheses). Theorem~\ref{th1} below shows that the almost sure convergence condition
\eqref{asgeneral} is sufficient (but not necessary) for obtaining lower bounds for all positive
moments of the detection delay. In fact, the condition \eqref{inprobk} in Lemma~\ref{le2} and
Theorem~\ref{th1} holds whenever $Z_{k+n-1}^k/n$ converges almost surely to the number $q$.

The following lemma will be used to derive asymptotic lower bounds for any positive moment of the
detection delay.


\begin{lemma} \label{le2} 
Let $Z_n^k$ be defined as in \eqref{LLRgen} and assume that for some $q>0$
\begin{equation} \label{inprobk}
\Pb_{k}\set{\frac{1}{M}\max_{1\le n \le M}Z_{k+n-1}^k\ge (1+\vae) q}\xra[M\to\infty]{}0 \quad
\text{for all $\vae>0$ and $k\ge 1$}.
\end{equation}
Then, for all $0<\vae<1$,
\begin{equation} \label{ga0general}
\lim_{\al\to 0} \sup_{\tau\in \Dal} \ga_{\vae,\al}^\pi(\tau)=0
\end{equation}
and
\begin{equation} \label{ga0tauA}
\lim_{A\to \infty} \ga_{\vae,A}^{\pi}(\tau_A)=0 .
\end{equation}

\end{lemma}

By \eqref{PktoPinfty},
\begin{equation} \label{strongerconstraint}
\Pb^\pi(\tau < \la)=\sum_{k=1}^\infty\pi_k \Pb_\infty(\tau<k) \le \Pb_\infty (\tau<\infty).
\end{equation}
Therefore, Lemma 1 of Tartakovsky and Veeravalli \cite{TartakovskyVeeravalli-TVP2004} may be
applied to prove statements \eqref{ga0general} and \eqref{ga0tauA} for the classes of prior
distributions considered in that work (i.e., for priors with exponential right tails and for
heavy-tailed priors). However, here we do not restrict ourselves to these classes of prior
distributions. The proof of the lemma for an arbitrary prior distribution is given in the Appendix.

Making use of Lemma~\ref{le2} and Chebyshev's inequality allows us to obtain the asymptotic lower
bounds for positive moments of the detection delay $\D_m^\pi(\tau)$, $m>0$.


\begin{theorem} \label{th1} 
Suppose condition \eqref{inprobk} holds for some positive finite number $q$. Then, for all $m >0$,
\begin{equation} \label{DmlowergentauA}
\D_m^\pi(\tau_A) \ge \brc{\frac{\log A}{q}}^m (1+o(1)) \quad \text{as $A\to\infty$}
\end{equation}
and
\begin{equation} \label{Dmlowergeninf}
\inf_{\tau\in\Dal}\D_m^\pi(\tau) \ge \brc{\frac{|\log \al|}{q}}^m (1+o(1)) \quad \text{as $\alo$},
\end{equation}
where $o(1)\to 0$.
\end{theorem}
%

\proof By the Chebyshev inequality, for any $0 < \vae <1$, $m>0$, and any $\tau \in \Dal$
\[
\Eb^\pi[(\tau-\la)^+]^m \ge[(1-\vae)L_\al]^m \Pb^\pi\set{\tau-\la \ge (1-\vae)L_\al},
\]
where
\[
\Pb^\pi\set{\tau-\la \ge (1-\vae)L_\al}=\Pb^\pi\set{\tau \ge \la}-\ga_{\vae,\al}^\pi(\tau).
\]
By \eqref{strongerconstraint}, for any $\tau\in\Dal$
\[
\Pb^\pi(\tau \ge \la) \ge 1- \Pb_\infty (\tau <\infty) \ge 1- \al.
\]
Thus, for any $\tau\in\Dal$
\begin{equation} \label{Dminflower}
\begin{split}
\D_m^\pi(\tau) & =\frac{\Eb^\pi[(\tau-\la)^+]^m}{\Pb^\pi\set{\tau\ge\la}} \ge[(1-\vae)L_\al]^m
\left[1-\frac{\ga_{\vae,\al}^\pi(\tau)}{\Pb^\pi\set{\tau\ge\la}}\right]
\\
& \ge [(1-\vae)L_\al]^m \left[1-\frac{\ga_{\vae,\al}^\pi(\tau)}{1-\al}\right].
\end{split}
\end{equation}
Since $\vae$ can be arbitrarily small and, by Lemma~\ref{le2},
$\sup_{\tau\in\Dal}\ga_{\vae,\al}^\pi(\tau)\to 0$ as $\al\to 0$, the asymptotic lower bound
\eqref{Dmlowergeninf} follows.

To prove \eqref{DmlowergentauA}, it suffices to repeat the above argument replacing $\al$ with
$1/A$ and using the fact that $\Pb_\infty(\tau_A < \infty) \le 1/A$ by Lemma \ref{le1}.
\endproof
%

Consider now the traditional i.i.d.\ model \eqref{eq:pdfiid} with pre-change and post-change
densities $f_0(x)$ and $f_1(x)$ (with respect to a sigma-finite measure $\mu(x)$), in which case
the LLR \eqref{LLRgen} is given by
\begin{equation} \label{LLRiid}
Z_n^{k} = \sum_{i=k}^n \log \frac{f_{1}(X_i)}{f_{0}(X_i)}, \quad k \le n.
\end{equation}
Define the Kullback-Leibler information number
\[
I = I(f_1,f_0)=\int \log \brc{\frac{f_1(x)}{f_0(x)}} f_1(x) \, d \mu(x),
\]
and assume that $0< I < \infty$. Then $\Eb_k Z_{k+n-1}^k= I \, n$ and the almost sure convergence
condition \eqref{asgeneral} holds with $q=I$ by the strong law of large numbers, i.e.,
\begin{equation}\label{asiid}
\frac{1}{n}Z_{k+n-1}^{k} \xra[n\to\infty]{\Pb_k-\text{a.s.}} I \quad \text{for every $k<\infty$}.
\end{equation}

Note that in the i.i.d.\ case condition \eqref{inprobk} holds with $q=I$. Therefore, as the first
step we have the following corollary that establishes the lower bound in the i.i.d.\ case.


\begin{corollary} \label{cor1} 
Let the Kullback-Leibler information number be positive and finite, $0<I<\infty$. Then the
asymptotic lower bounds \eqref{DmlowergentauA} and \eqref{Dmlowergeninf} hold with $q=I$.
\end{corollary}

\subsection{Asymptotic optimality in the i.i.d.\ case} \label{ss:AOiid}

We now proceed with devising first-order approximations to the moments of the detection delay of
the detection test $\tau_A$ as $A\to\infty$ and establishing its first-order asymptotic optimality
when $A=1/\al$ and $\al\to 0$ in the i.i.d.\ case.

\medskip

\subsubsection{First-order approximations} \label{sss:FO}

In order to prove the asymptotic optimality property, it suffices to derive an upper bound showing
that this bound is asymptotically the same as the lower bound specified in Corollary \ref{cor1}.

It is easily seen that for any $k \ge 1$
\begin{align}
G_n & = \Pi_{n+1} + \sum_{j=1}^n \pi_j e^{Z_n^j} \nonumber
\\
& = e^{Z_n^k} \brc{\pi_k + \Pi_{n+1} e^{-Z_n^k}+ \sum_{j=1}^{k-1} \pi_j e^{\sum_{i=j}^{k-1} \De
Z_i} + \sum_{j=k}^{n-1} \pi_{j+1} e^{-\sum_{i=k}^{j} \De Z_i}} \label{eq:GnZnk}
\\
& \ge e^{Z_n^k} \pi_k , \label{ineq:GnZnk}
\end{align}
where $\De Z_i= \log [f_1(X_i)/f_0(X_i)]$. Thus, for any $k \ge 1$, the stopping
time $\tau_A$ does not exceed the stopping time
\begin{equation} \label{nuk}
\nu_k(A) = \min \set{n \ge k: Z_n^k \ge \log(A/\pi_k)}.
\end{equation}
Moreover,
\[
(\tau_A-k)^+ \le \nu_k(A)-k .
\]
By the i.i.d.\ property of the data, the random variables $\De Z_i$, $i =1,2,\dots$
are also i.i.d.\ and hence the distribution of $\nu_k(A)-k+1$ under $\Pb_k$ is the
same as the $\Pb_1$-distribution of the stopping time
\begin{equation} \label{tildenuk}
\tilde{\nu}_1(A,\pi_k) = \min \set{n \ge 1: Z_n^1 \ge \log (A/\pi_k)}.
\end{equation}
Therefore, for all $k \ge 1$
\begin{equation} \label{EtauAleEnuk}
\Eb_k[(\tau_A-k)^+]^m  \le \Eb_1 (\tilde{\nu}_1(A,\pi_k)-1)^m ,
\end{equation}
which can be used to obtain the desired upper bound.

Details are given in the following theorem.

\begin{theorem} \label{th2} 
Let $0<I <\infty$ and let prior distribution be such that $\sum_{k=1}^\infty |\log \pi_k|^m \pi_k <
\infty$.

(i) As $A \to \infty$,
\begin{equation} \label{DmtauAiid}
\D_m^\pi(\tau_A) \sim \brc{\frac{\log A}{I}}^m.
\end{equation}

(ii) If $A_\al=1/\al$, then $\tau_{A_\al}\in \Dal$ and, as $\al\to 0$, for all $m
\ge 1$
\begin{equation} \label{DmAOiid}
\inf_{\tau\in\Dal}\D_m^\pi(\tau) \sim \D_m^\pi(\tau_{A_\al}) \sim \brc{\frac{|\log
\al|}{I}}^m .
\end{equation}
\end{theorem}
%


\proof
(i) In the i.i.d.\ case, the LLR $Z_{n}^1$, $n \ge 1$ is a random walk with mean
$\Eb_1 Z_n^1= I\, n$. Since $I$ is positive and finite, $\Eb_1
\{-\min(0,Z_1^1)\}^m<\infty$ for all $m>0$. Indeed,
\[
\Eb_1 \exp\set{-\min(0,Z_1^1)}= \Eb_1 e^{-Z_1^1}\Ind{Z_1^1<0} + \Eb_1 \Ind{Z_1^1 \ge
0} \le \Eb_1 e^{-Z_1^1} +1 = 2 .
\]
Therefore, we can apply Theorem~III.8.1 of Gut~\cite{Gut} that yields, for all $m\ge
1$,
\begin{equation} \label{Enuk}
\Eb_1 [\tilde{\nu}_1(A,\pi_k)]^m = \brc{\frac{\log(A/\pi_k)}{I}}^m (1+o(1)) \quad
\text{as $A \to\infty$}.
\end{equation}
Using \eqref{Enuk} along with \eqref{EtauAleEnuk} implies
\begin{equation} \label{EtauAkupper}
\Eb_k[(\tau_A-k)^+]^m \le \brc{\frac{\log (A/\pi_k)}{I}}^m [1+ \vae(k,m,A)]
\end{equation}
where $\vae(k,m,A) \to 0$ as $A\to\infty$.

Write $a=\log A$. Now, averaging in \eqref{EtauAkupper} over the prior distribution, we obtain
\begin{equation} \label{Averageineq}
\begin{split}
\sum_{k=1}^\infty \pi_k \Eb_k[(\tau_A-k)^+]^m & \le \brc{\frac{a}{I}}^m \Bl\{\sum_{k=1}^\infty
\pi_k \Bl(1+ \frac{|\log \pi_k|}{a}\Br)^m
\\
& \quad + \sum_{k=1}^\infty \pi_k \Bl(1+ \frac{|\log \pi_k|}{a}\Br)^m \vae(k,m,A)\Br\} \quad
\text{as $A \to\infty$},
\end{split}
\end{equation}
Since by the conditions of the theorem $\sum_{k=1}^\infty |\log \pi_k|^m \pi_k < \infty$, it
follows that
\begin{equation} \label{firstto0}
\sum_{k=1}^\infty \pi_k \brc{1+ \frac{|\log \pi_k|}{a}}^m =1+o(1) \quad \text{as $A\to\infty$}.
\end{equation}
The important observation is that since $|\log \pi_k| \to \infty$ as $k\to\infty$, the asymptotic
equality \eqref{Enuk} and, hence, the inequality \eqref{EtauAkupper} also hold for any $A>1$ as $k
\to\infty$. This means that $\vae(k,m,A) \to 0$ as $k\to\infty$ for any fixed $A>1$ and also as
$A\to\infty$. It follows that
\[
\sum_{k=1}^\infty \pi_k |\log \pi_k|^m \vae(k,m,A)<\infty \quad \text{for any $A>1$}
\]
and, hence,
\begin{equation} \label{secondto0}
\sum_{k=1}^\infty \pi_k \brc{1+ \frac{|\log \pi_k|}{a}}^m \vae(k,m,A) \to 0 \quad \text{as
$A\to\infty$}.
\end{equation}

Combining \eqref{Averageineq}, \eqref{firstto0}, and \eqref{secondto0} yields the asymptotic
inequality
\[
\Eb^\pi[(\tau_A-\la)^+]^m \le \brc{\frac{\log A}{I}}^m (1+o(1)) \quad \text{as $A \to\infty$}.
\]

Finally, noting that $\Pb^\pi (\tau_A \ge \la) \ge \Pb_\infty (\tau_A =\infty) \ge
1- 1/A$ (cf.\ Lemma \ref{le1}) and
\[
\Eb^\pi[(\tau_A-\la)^+]^m = \Pb^\pi (\tau_A \ge \la) \D_m^\pi(\tau_A),
\]
we obtain the upper bound
\[
\D_m^\pi(\tau_A) \le \brc{\frac{\log A}{I}}^m (1+o(1)).
\]
Comparing this asymptotic upper bound with the lower bound \eqref{DmlowergentauA} (see
Corollary~\ref{cor1}) completes the proof of \eqref{DmtauAiid}.

(ii) The fact that $\tau_{A_\al} \in \Dal$ when $A_\al=1/\al$ follows from Lemma
\ref{le1}. The asymptotic relation \eqref{DmAOiid} follows from \eqref{DmtauAiid}
and the lower bound \eqref{Dmlowergeninf}.

\endproof


\subsubsection{Higher-order approximations} \label{sss:HO}
The upper bound $\Pb_\infty(\tau_A < \infty) \le 1/A$ (see \eqref{boundPFA}) for the global PFA,
which neglects a threshold overshoot, holds in the most general, non-i.i.d.\ case. In the i.i.d.\
case, an accurate approximation for $\Pb_\infty(\tau_A < \infty)$ can be obtained by taking into
account an overshoot using the nonlinear renewal theory argument (see
Woodroofe~\cite{Woodroofebook82} and Siegmund~\cite{Siegmundbook85}). This is important in
situations where the upper bound \eqref{boundPFA} that ignores the overshoot is conservative, which
is always the case where the densities $f_1(x)$ and $f_0(x)$ are not close enough.

In order to apply relevant results from nonlinear renewal theory, we have to rewrite the stopping
time $\tau_A$ in the form of a random walk crossing a constant threshold plus a nonlinear term that
is slowly changing in the sense defined in \cite{Siegmundbook85,Woodroofebook82}. Using
\eqref{eq:GnZnk} and writing
\begin{equation} \label{ell}
\ell_n^k = \log \brc{\pi_k+ \Pi_{n+1} e^{-\sum_{i=k}^n \De Z_i} + \sum_{j=1}^{k-1} \pi_j
e^{\sum_{i=j}^{k-1} \De Z_i} + \sum_{j=k}^{n-1} \pi_{j+1} e^{-\sum_{i=k}^j \De Z_i}},
\end{equation}
we obtain that for every $k \ge 1$
\begin{equation} \label{logGn}
\log G_n = Z_n^k + \ell_n^k .
\end{equation}
Therefore, on $\{\tau_A \ge k\}$ for any $k \ge 1$, the stopping time $\tau_A$ can be written in
the following form
\begin{equation} \label{nusc}
\tau_A=\min\set{n\ge k: Z_n^k + \ell_n^k  \ge a}, \quad a=\log A,
\end{equation}
where $\ell_n^k$ is given by \eqref{ell} and $Z_n^k$, $n \ge k$ is a random walk with mean $\Eb_k
Z_n^k= I\, n$.

For $b>0$, define $\eta_{b}$ as
\begin{equation} \label{eta1}
\eta_{b}= \min\{n\ge 1: Z_n^1 \ge b\} ,
\end{equation}
and let $\varkappa_b= Z_{\eta_b}^1-b$ (on $\{\eta_b<\infty\}$) denote the excess
(overshoot) of the statistic $Z_n^1$ over the threshold $b$ at time $n=\eta_b$. Let
\begin{equation} \label{limdisover}
H(y,I)=\lim_{b\to\infty}\Pb_1\set{\varkappa_b \le y}
\end{equation}
be the limiting distribution of the overshoot and let
\begin{equation} \label{averageexpovershoot}
\zeta(I)=\lim_{b\to\infty}\Eb_1 e^{-\varkappa_b}=\int_0^\infty e^{-y} \; d H(y,I).
\end{equation}

The important observation is that $\ell_n^k$, $n\ge 1$ are slowly changing. To see this it suffices
to note that, as $n\to\infty$, the values of $\ell_n^k$ converge to the random variable
\[
\ell_\infty^k=\log \brc{\pi_k + \sum_{j=1}^{k-1} \pi_j e^{\sum_{i=j}^{k-1} \De Z_i} +
\sum_{j=k}^{\infty} \pi_{j+1} e^{-\sum_{i=k}^j \De Z_i}},
\]
which has finite negative expectation. Indeed, on the one hand $\ell_\infty^k \ge \log \pi_k$, and
on the other hand, by Jensen's inequality,
\begin{align*}
\Eb_k \ell_\infty^k & \le \log \brc{\pi_k + \sum_{j=1}^{k-1} \pi_j \Eb_k e^{\sum_{i=j}^{k-1} \De
Z_i} + \sum_{j=k}^{\infty} \pi_{j+1} \Eb_k e^{-\sum_{i=k}^j \De Z_i}}
\\
& = \log \brc{\pi_k + \sum_{j=1}^{k-1} \pi_j  + \sum_{j=k}^{\infty} \pi_{j+1}} =\log
\brc{\sum_{j=1}^{\infty} \pi_{j}} =0,
\end{align*}
where we used the equalities
\[
\Eb_k e^{\sum_{i=j}^{k-1} \De Z_i}= \prod_{i=j}^{k-1} \Eb_k \frac{f_1(X_i)}{f_0(X_i)} =1
\]
and
\[
\Eb_k e^{-\sum_{i=k}^j \De Z_i} = \prod_{i=k}^{j} \Eb_k \frac{f_0(X_i)}{f_1(X_i)} =1,
\]
which hold since, obviously,
\[
\Eb_k \frac{f_1(X_i)}{f_0(X_i)} = \int  \frac{f_1(x)}{f_0(x)} f_0(x) d\mu(x) =1 \quad \text{for $i
<k$}
\]
and
\[
\Eb_k \frac{f_0(X_i)}{f_1(X_i)} = \int  \frac{f_0(x)}{f_1(x)} f_1(x) d\mu(x) =1 \quad \text{for $i
\ge k$}.
\]

An important consequence of the slowly changing property is that, under mild conditions, the
limiting distribution of the overshoot of a random walk does not change by the addition of a slowly
changing nonlinear term (see Theorem~4.1 of Woodroofe \cite{Woodroofebook82}). This property allows
us to derive an accurate asymptotic approximation for the probability of false alarm, which is
important in situations where the value of $I$ is moderate. (For small values of $I$ the overshoot
can be neglected.) The following lemma presents an exact result.

%

\begin{lemma} \label{le3} 
Suppose $Z_n^1$, $n\ge 1$ are nonarithmetic with respect to $\Pb_1$. Let $I <\infty$. Then
\begin{equation} \label{HOPFA}
\Pb_\infty(\tau_A < \infty)=\frac{\zeta(I)}{A} (1+o(1)) \quad \text{as
$A\to\infty$}.
\end{equation}
\end{lemma}
%


\proof

Obviously,
\begin{align*}
\Pb_\infty(\tau_A < \infty) & =\Eb^\pi\set{G_{\tau_A}^{-1} \Ind{\tau_A < \infty}}=
\Eb^\pi(\set{A/(A G_{\tau_A}) \Ind{\tau_A < \infty}}=
\\
&  = \frac{1}{A}\Eb^\pi \set{e^{-\chi_a}\Ind{\tau_A < \infty}},
\end{align*}
where $\chi_a=\log G_{\tau_A}-a$. Since $\chi_a \ge 0$ and $\Pb^\pi(\tau_A < \la)
\le \Pb_\infty(\tau_A < \infty)\le 1/A$, it follows that
\begin{align*}
\Eb^\pi \set{e^{-\chi_a}\Ind{\tau_A < \infty}} &
=\Eb^\pi\set{e^{-\chi_a}|\tau_A<\la} \Pb^\pi(\tau_A < \la) +
\Eb^\pi\set{e^{-\chi_a}|\tau_A\ge\la} \brc{1-\Pb^\pi(\tau_A < \la)}
\\
& = \Eb^\pi\set{e^{-\chi_a}|\tau_A\ge\la} +O(1/A) \quad \text{as $A\to\infty$}.
\end{align*}

Therefore, it suffices to evaluate the value of
\[
\Eb^\pi\set{e^{-\chi_a}|\tau_A\ge\la} = \sum_{k=1}^\infty \Pb(\la=k| \tau_A \ge k)
\Eb_k\set{e^{-\chi_a}|\tau_A\ge k}.
\]
Recall that, by \eqref{nusc}, for any $1 \le k <\infty$,
\begin{equation*} 
\tau_A=\set{n\ge k: Z_n^k+\ell_n^k \ge a} \quad \text{on $\{\tau_A \ge k\}$},
\end{equation*}
where $Z_n^k$, $n\ge k$ is a random walk with the expectation $\Eb_k Z_k^k = I$ and $\ell_n^k$, $n
\ge k$ are slowly changing under $\Pb_k$. Since, by conditions of the lemma, $0< I<\infty$, we can
apply Theorem~4.1 of Woodroofe~\cite{Woodroofebook82} to obtain
\[
\lim_{A\to\infty} \Eb_k\set{e^{-\chi_a}|\tau_A\ge k}=\int_0^\infty e^{-y} \; dH(y, I)= \zeta(I).
\]
Since $\Pb_\infty (\tau_A \ge k) \ge 1-1/A$ and $\Pb^\pi (\tau_A \ge \la) \ge
1-1/A$,
\[
\lim_{A\to\infty}\Pb(\la=k| \tau_A \ge k) = \lim_{A\to\infty} \frac{\pi_k \Pb_\infty
(\tau_A \ge k)}{\Pb^\pi (\tau_A \ge \la)} =\pi_k
\]
and, therefore,
\[
\lim_{A\to\infty} \Eb^\pi\set{e^{-\chi_a}|\tau_A\ge \la}=\lim_{A\to\infty}
\Eb^\pi\set{e^{-\chi_a}}=\zeta(I),
\]
which completes the proof of \eqref{HOPFA}.

\endproof


Under an additional, second moment condition, the nonlinear renewal theorem
\cite{Woodroofebook82} also allows for obtaining a higher-order approximation for
the ADD:
\begin{align}
\Eb_k(\tau_A-k|\tau_A\ge k) & = I^{-1} \bl[\log A - C_k^\pi(I) +
\overline{\varkappa}(I)\br] + o(1), \quad k \ge 1; \label{CADDktauA}
\\
\ADD^{\pi}(\tau_A) = \Eb^\pi(\tau_A-\la|\tau_A\ge \la) & = I^{-1} \bl[\log A - \sum_{k=1}^\infty
C_k^\pi(I) \pi_k + \overline{\varkappa}(I)\br] + o(1), \label{ADDtauA}
\end{align}
where $C_k^\pi(I)=\Eb_k\ell_\infty^k$ and
\begin{equation} \label{averageovershoot}
\overline{\varkappa}(I)=\lim_{a\to\infty}\Eb_1 \varkappa_a=\int_0^\infty y \; d H(y,I)
\end{equation}
is the limiting average overshoot in the one-sided test.

However, approximations \eqref{CADDktauA} and \eqref{ADDtauA} have little value,
since it is usually impossible to compute the constant $C_k^\pi(I)$. Instead, we
propose the following approximations
\begin{align}
\Eb_k(\tau_A-k|\tau_A\ge k) & \approx I^{-1} \bl[\log (A/\pi_k) + \overline{\varkappa}(I)-1 \br],
\quad k \ge 1;\label{HOCADDktauAappr}
\\
\ADD^{\pi}(\tau_A) & \approx I^{-1} \bl[\log A + \sum_{k=1}^\infty \pi_k |\log \pi_k| +
\overline{\varkappa}(I)-1\br] , \label{HOADDtauAappr}
\end{align}
which use the minimal value of the random variable $\ell_\infty^k=\log \pi_k$. Clearly, one may
expect that these approximations will overestimate the true values. On the other hand, it is
expected that the approximations that ignore the overshoot given by
\begin{align}
\Eb_k(\tau_A-k|\tau_A\ge k) & \approx I^{-1} \bl[\log (A/\pi_k) -1 \br] , \quad k \ge 1;
\label{FOCADDktauAappr}
\\
\ADD^{\pi}(\tau_A) & \approx I^{-1} \bl[\log A + \sum_{k=1}^\infty \pi_k |\log \pi_k| -1 \bl]
\label{FOADDtauAappr}
\end{align}
will underestimate the true values.

The constants $\zeta(I)$ and $\overline{\varkappa}(I)$ defined in
\eqref{averageexpovershoot} and \eqref{averageovershoot} are the subject of the
renewal theory. They can be computed either exactly or approximately in a variety of
particular examples.

\subsection{Asymptotic optimality in the non-i.i.d.\ case} \label{ss:AOnoniid}

In this section, we deal with the general non-i.i.d.\ model \eqref{eq:pdfnoniid} and show that
under certain quite general conditions the detection procedure \eqref{eq:tauA} is asymptotically
optimal for small $\al$.

As we established in Theorem~\ref{th1} above, the strong law of large numbers \eqref{asgeneral} is
sufficient for obtaining the lower bound for the moments of the detection delay. However, in
general, this condition is not sufficient for asymptotic optimality with respect to the moments of
the detection delay.  Therefore, some additional conditions are needed to guarantee asymptotic
optimality.

\medskip

\subsubsection{Weak asymptotic optimality} \label{sss:WeakAO}

We begin with answering the question of whether some asymptotic optimality result can still be
obtained under the almost sure convergence condition \eqref{asgeneral}. The following theorem
establishes asymptotic optimality of the procedure $\tau_{A_\al}$ in a weak probabilistic sense.


\begin{theorem} \textbf{(Weak Asymptotic Optimality)} \label{th3} 
Let there exist a finite positive number $q$ such that condition \eqref{asgeneral} hold, and let
$A=A_\al=1/\al$. Then, for every $0<\vae<1$,
\begin{equation} \label{weakoptk}
\inf_{\tau\in\Dal}\Pb_k\set{(\tau-k)^+ >\vae (\tau_{A_\al}-k)^+}\xra[\al\to 0]{} 1 \quad \text{for
all $k \ge 1$}
\end{equation}
and
\begin{equation}  \label{weakoptaverage}
\inf_{\tau\in\Dal} \Pb^\pi\set{(\tau-\la)^+ >\vae (\tau_{A_\al}-\la)^+}\xra[\al\to 0]{} 1.
\end{equation}
\end{theorem}


\proof

Extracting the term $e^{Z_n^k}$, the statistic $G_n$ can be written as follows:
\begin{equation} \label{Gngen}
G_n = \pi_k e^{Z_n^k} \brc{1 + \frac{1}{\pi_k} \left[ \Pi_{n+1} e^{-Z_n^k} + \sum_{j=1}^{k-1} \pi_j
e^{Z_{k-1}^j} + \sum_{j=k}^{n-1} \pi_{j+1} e^{-Z_{k}^j}\right ]} .
\end{equation}
Writing
\[
Y_n^k = \pi_k^{-1} \brc{\Pi_{n+1} e^{-Z_n^k} + \sum_{j=1}^{k-1} \pi_j e^{Z_{k-1}^j} +
\sum_{j=k}^{n-1} \pi_{j+1} e^{-Z_{k}^j}},
\]
we obtain that for every $k \ge 1$
\begin{equation} \label{logGnnew}
\log G_n = Z_n^k + \log (1+ Y_n^k) + \log \pi_k.
\end{equation}
It is easily verified that $\Eb_k e^{-Z_n^k}=1$, $\Eb_k e^{-Z_k^j}=1$ for $j \ge k$, and $\Eb_k
e^{Z_{k-1}^j}=1$ for $j \le k-1$ and, hence,
\[
\Eb_k Y_n^k = \pi_k^{-1} \brc{\Pi_{n+1} + \sum_{j=1}^{k-1} \pi_j + \sum_{j=k}^{n-1} \pi_{j+1}} =
(1-\pi_k)/\pi_k.
\]

Since $\log (1+ Y_n^k)$ is non-negative, applying Markov's inequality we obtain that for every
$\vae >0$
\[
\Pb_k\set{n^{-1}\log (1+Y_n^k) \ge \vae} \le e^{-n\vae}(1+ \Eb_k Y_n^k) =e^{-n\vae}/
\pi_k.
\]
It follows that for all $\vae>0$
\[
\sum_{n=k}^\infty \Pb_k\set{n^{-1}\log (1+Y_n^k) \ge \vae} < \infty,
\]
which implies that
\begin{equation} \label{logYas}
\frac{1}{n-k+1} \log (1+ Y_n^k) \xra[n\to\infty]{\Pb_k-\text{a.s.}} 0 \quad \text{for every $k \ge
1$}.
\end{equation}
Using \eqref{asgeneral}, \eqref{logGnnew}, and \eqref{logYas} yields
\begin{equation} \label{logGas}
\frac{1}{n} \log G_n \xra[n\to\infty]{\Pb_k-\text{a.s.}} q \quad \text{for every $k
\ge 1$}.
\end{equation}

Clearly, $\tau_A \to \infty$ as $A\to \infty$ almost surely under $\Pb_k$ for every
$k \ge 1$ and, by \eqref{logGas}, $G_n\to \infty$ a.s.\ under $\Pb_k$ (as $n \to
\infty$), which implies that $\Pb_k(\tau_A < \infty)=1$. Therefore,
\begin{equation} \label{ineqlogBovertauA}
q \xla[A\to\infty]{\Pb_k-\text{a.s.}} \frac{\log G_{\tau_A-1}}{\tau_A} \le
\frac{\log A}{\tau_A} \le \frac{\log G_{\tau_A}}{\tau_A}
\xra[A\to\infty]{\Pb_k-\text{a.s.}} q
\end{equation}
and, since $\Pb_k (\tau_A < k) \le 1/A \to 0$, it follows that
\begin{equation} \label{tauAinprobk}
\frac{(\tau_{A_\al}-k)^+}{|\log \al|} \to \frac{1}{q} \quad \text{in
$\Pb_k-$probability as $\al\to 0$ for all $k\ge 1$}
\end{equation}
and
\begin{equation} \label{tauAinprobaverage}
\frac{(\tau_{A_\al}-\la)^+}{|\log \al|} \to \frac{1}{q} \quad \text{in
$\Pb^\pi-$probability as $\al\to 0$.}
\end{equation}
Next, since the right side in inequality \eqref{gakupper} (see Appendix) does not depend on the
stopping time $\tau$ it follows
\begin{equation} \label{probkinf}
\lim_{\al\to 0} \inf_{\tau\in \Dal} \Pb_k \set{\tau -k \ge \vae q^{-1}|\log \al|} =
1 \quad \text{for all $k \ge 1$ and $0<\vae <1$},
\end{equation}
which along with \eqref{tauAinprobk} proves \eqref{weakoptk}.

Finally, the asymptotic relation \eqref{weakoptaverage} follows from
\eqref{tauAinprobaverage} and Lemma~\ref{le2}, which implies that
\[
\lim_{\al\to 0} \inf_{\tau\in \Dal} \Pb^\pi \set{\tau -\la \ge \vae q^{-1}|\log
\al|} = 1 \quad \text{for all $0<\vae <1$}.
\]

\endproof


\subsubsection{First-order asymptotic optimality} \label{sss:FOAOnoniid}

We now proceed with the first-order (FO) asymptotic optimality with respect to positive moments of
the detection delay $\D_m^\pi(\tau)$. We first note that using the method proposed by
Lai~\cite{LaiIEEE98} it can be shown that the ADD of the detection procedure $\tau_{A_\al}$ attains
the lower bound \eqref{Dmlowergeninf} ($m=1$) under the condition
\[
\max_{1\le k \le j} \Pb_k \set{n^{-1} Z_{j+n-1}^j \le q-\vae} \to 0 \quad \text{as
$n\to\infty$ for all $\vae >0$}.
\]
It can be also shown that, for any $m \le r$, the sufficient condition for
$\D_m^\pi(\tau_{A_\al})$ to attain the lower bound \eqref{Dmlowergeninf} is
\[
\sum_{k=1}^\infty \sum_{n=1}^\infty \pi_k n^{r-1} \Pb_k \set{Z_{k+n-1}^k \le
(q-\vae)n} < \infty \quad \text{for all $\vae > 0$}.
\]
This latter condition is closely related to the following condition
\begin{equation} \label{rquicklefttail}
\sum_{k=1}^\infty \pi_k \Eb_k (T_{k,\vae})^r < \infty \quad \text{for all $\vae >
0$},
\end{equation}
where
\begin{equation} \label{Tk}
T_{k,\vae}=\sup\set{n\ge 1:n^{-1} Z_{k+n-1}^k - q < -\vae} \quad (\sup\set{\vano}=0)
\end{equation}
is the last time when $n^{-1} Z_{k+n-1}^k$ leaves the region $[q-\vae, \infty)$.


\begin{theorem} \textbf{(FO Asymptotic Optimality)}  \label{th4} 
Let conditions \eqref{inprobk} and \eqref{rquicklefttail} hold for some positive finite $q$ and
some $r \ge 1$. Assume that
\begin{equation} \label{prior}
\sum_{k=1}^\infty |\log \pi_k|^m \pi_k < \infty \quad \text{for $m \le r$}.
\end{equation}
Then for all $m \le r$
\begin{equation} \label{mthmomenttau}
\D_m^\pi(\tau_A) \sim \brc{\frac{|\log A|}{q}}^m \quad \text{as $A \to 0$}.
\end{equation}

If $A=A_\al=1/\al$, then $\tau_{A_\al}\in \Dal$ and for all $m \le r$,
\begin{equation} \label{mthmomentAO}
\inf_{\tau \in \Dal} \D_m^\pi(\tau) \sim \D_m^\pi(\tau_{A_\al}) \sim
\brc{\frac{|\log\al|}{q}}^m \quad \text{as $\al\to 0$}.
\end{equation}
\end{theorem}
%


\proof

To prove \eqref{mthmomenttau} it suffices to show that the lower bound \eqref{DmlowergentauA} in
Theorem~\ref{th1} is also asymptotically the upper bound, i.e.,
\begin{equation} \label{DmuppergentauA}
\D_m^\pi(\tau_A) \le \brc{\frac{\log A}{q}}^m (1+o(1)) \quad \text{as $A\to\infty$}.
\end{equation}

It follows from equality \eqref{Gngen} that
\[
G_n  \ge e^{Z_n^k} \pi_k ,
\]
and, therefore, for any $k \ge 1$,
\begin{equation} \label{taulenu}
(\tau_A -k)^+ \le \nu_k(A) = \min \set{n \ge 1: Z_{k +n-1}^k \ge \log (A/\pi_k)}.
\end{equation}
Thus,
\[
\D_m^\pi(\tau_A) \le \frac{\sum_{k=1}^\infty \pi_k \Eb_k (\nu_k(A))^m}{\Pb (\tau_A
\ge \la)} .
\]
Since by Lemma \ref{le1} $\Pb (\tau_A \ge \la) \ge 1- \Pb_\infty (\tau_A < \infty) \ge 1-1/A$, it
is sufficient to prove that
\begin{equation} \label{Emnuupper}
\sum_{k=1}^\infty \pi_k \Eb_k (\nu_k(A))^m \le \brc{\frac{\log A}{q}}^m (1+o(1))
\quad \text{as $A\to\infty$}.
\end{equation}

By the definition of the stopping time $\nu_k$,
\[
Z_{k+\nu_k-2}^k < \log (A/\pi_k) \quad \text{on $\set{\nu_k < \infty}$}.
\]
On the other hand, by the definition of the last entry time \eqref{Tk},
\[
Z_{k+\nu_k-2}^k \ge (q-\vae)(\nu_k-1) \quad \text{on $\set{\nu_k > 1+ T_{k,\vae}}$}.
\]
Hence,
\[
(q-\vae)(\nu_k-1) \le \log (A/\pi_k) \quad \text{on $\set{T_{k,\vae} +1 < \nu_k < \infty}$}
\]
and we obtain
\[
\Eb_k \nu_k^m = \Eb_k \nu_k^m \Ind{T_{k,\vae} +1 < \nu_k < \infty} + \Eb_k \nu_k^m
\Ind{\nu_k \le T_{k,\vae} +1} \le \brc{\frac{1+ \log (A/\pi_k)}{q-\vae}}^m + \Eb_k
(1+T_{k,\vae})^m.
\]
Averaging over the prior distribution yields
\[
\sum_{k=1}^\infty \pi_k \Eb_k \nu_k^m \le \sum_{k=1}^\infty \pi_k \brc{\frac{1+ \log
(A/\pi_k)}{q-\vae}}^m + \sum_{k=1}^\infty \pi_k \Eb_k (1+T_{k,\vae})^m .
\]
By conditions \eqref{prior} and \eqref{rquicklefttail}, $\sum_{k=1}^\infty |\log
\pi_k|^m\pi_k<\infty$ and $\sum_{k=1}^\infty \pi_k \Eb_k (T_{k,\vae})^m < \infty$ for $m \le r$.
Since $\vae$ can be arbitrarily small, the asymptotic upper bound \eqref{Emnuupper} follows and the
proof of \eqref{mthmomenttau} is complete.

Asymptotic relations \eqref{mthmomentAO} follow from \eqref{mthmomenttau} and the
asymptotic lower bound \eqref{Dmlowergeninf} in Theorem~\ref{th1}.

\endproof


Introduce now the double-sided last entry time
\begin{equation} \label{Tdsk}
T_{k,\vae}^{ds}=\sup\set{n\ge 1: |n^{-1}Z_{k+n-1}^k - q| > \vae} \quad
(\sup\set{\vano}=0),
\end{equation}
which is the last time when $n^{-1} Z_{k+n-1}^k$ leaves the region $[q-\vae,
q+\vae]$. In terms of $T_{k,\vae}^{ds}$, the almost sure convergence of
\eqref{asgeneral} may be written as $\Pb_k \{T_{k,\vae}^{ds}<\infty\}=1$ for all
$\vae>0$ and $k\ge 1$, which implies condition \eqref{inprobk}.

If instead of condition \eqref{rquicklefttail} we impose the condition
\begin{equation} \label{rquickds}
\sum_{k=1}^\infty \pi_k \Eb_k (T_{k,\vae}^{ds})^r < \infty \quad \text{for all $\vae
> 0$ and some $r \ge 1$}
\end{equation}
that limits the behavior of both tails of the distribution of the LLR $Z_{k+n-1}^k$,
then both conditions \eqref{inprobk} and \eqref{rquicklefttail} are satisfied and,
therefore, the following corollary holds.


\begin{corollary} \label{cor2} 
Suppose condition \eqref{rquickds} is satisfied for some positive finite $q$. Then
the asymptotic relations \eqref{mthmomenttau} and \eqref{mthmomentAO} hold.
\end{corollary}

Note that the condition $\Eb_k (T_{k,\vae}^{ds})^r <\infty$ is not more than the so-called
$r$-quick convergence of $n^{-1}Z_{k+n-1}^k$ to $q$ under $\Pb_k$ (cf.\
Lai~\cite{Lairquick,LaiAS1981} and Tartakovsky \cite{TarSISP98}). It is closely related to the
condition
\[
\sum_{n=1}^\infty n^{r-1} \Pb_k\Bigl\{\Bigl|Z_{k+n-1}^k-q n\Bigr|>\vae n \Bigr\}<\infty \quad
\text{for all $\vae>0$},
\]
which determines the rate of convergence in the strong law of large numbers (cf.\ Baum and Katz
\cite{BaumKatz65} in the i.i.d.\ case). For $r=1$, the latter condition is the complete convergence
of $n^{-1}Z_{k+n-1}^{k}$ to $q$ under $\Pb_k$ (cf.\ Hsu and Robbins~\cite{HsuRobbins47}).

In particular examples, instead of checking the original condition \eqref{rquickds}, one may check
the following condition
\begin{equation} \label{rquickdsprob}
\sum_{n=1}^\infty \sum_{k=1}^\infty n^{r-1} \pi_k \Pb_k\Bigl\{\Bigl|Z_{k+n-1}^k-q n\Bigr|>\vae n
\Bigr\}<\infty \quad \text{for all $\vae>0$},
\end{equation}
which is sufficient for the asymptotic optimality property.

\begin{remark} \label{rem1}
In the i.i.d.\ case, the finiteness of the $(r+1)$-st absolute moment of the LLR,
$\Eb_1|Z_1^1|^{r+1} < \infty$, is both necessary and sufficient condition for the $r$-quick
convergence \eqref{rquickds}. See, e.g., Baum and Katz \cite{BaumKatz65}. Therefore,
Theorem~\ref{th4} implies asymptotic relations \eqref{DmtauAiid} and \eqref{DmAOiid} for $m \le r$
under the $(r+1)$-st moment condition. On the other hand, Theorem~\ref{th2} shows that these
relations hold for all $m >0$ under the unique first moment condition: $I < \infty$.
\end{remark}

\begin{remark} \label{rem2}
The asymptotic approximation \eqref{mthmomenttau} for the ADD ($m=1$) ignores the constant
$C_\pi=\sum_{k=1}^\infty\pi_k |\log \pi_k|$. The proof suggests that preserving this constant may
improve the accuracy of the first-order approximation for the ADD, i.e., the following approximate
formula
\[
\ADD^{\pi}(\tau_A) \approx q^{-1} (\log A + C_\pi)
\]
may be more accurate in particular examples.
\end{remark}

\section{Examples} \label{s:Examples}

\subsection{Detection of a change in the i.i.d.\ exponential sequence} \label{ss:Exp}

Let, conditioned on $\la=k$, the observations $X_1,\dots,X_{k-1}$ are i.i.d.\ $\Exp(1)$ and
$X_{k},X_{k+1},\dots$ are i.i.d.\ $\Exp(1/(1+Q))$, i.e.,
\[
f_1(x)= \frac{1}{1+Q} e^{-x/(1+Q)} \Ind{x\ge 0}, \quad f_0(x) = e^{-x} \Ind{x\ge 0},
\]
where  $Q >0$. Then the partial LLR $\De Z_n = - \log (1+Q) +[Q/(1+Q)] X_n$ and the
Kullback-Leibler information number
\[
I = \log (1+Q) - Q/(1+Q) .
\]
By Theorem \ref{th2}, the detection test $\tau_{A_\al}$ with $A_\al=1/\al$  minimizes
asymptotically as $\al\to 0$ all positive moments of the detection delay.

The distributions of the overshoot $\varkappa_b=Z^1_{\eta_b}-b$ in the one-sided, open-ended test
$\eta_b$ are exponential for all positive $b$ \cite{TartakovskyIvanovaPPI92}:
\[
\Pb_1 (\varkappa_b >x) = e^{-x/Q} \Ind{x \ge 0}, \quad \Pb_\infty (\varkappa_b >x) = e^{-x(1+Q)/Q}
\Ind{x \ge 0}
\]
and, therefore,
\begin{equation*}
\zeta(Q) = 1/(1+Q),  \quad \bar{\varkappa}(Q) = Q.
\end{equation*}
Note that these formulas are exact for any positive $b$, not just asymptotically as $b\to\infty$.

By Lemma \ref{le3}, if the threshold is set as
\[
A_\al = \frac{1}{(1+Q)\al},
\]
then for small $\al$
\[
\Pb_\infty(\tau_{A_\al}<\infty) = \al(1+o(1)),
\]
and by \eqref{HOADDtauAappr},
\[
\ADD^{\pi}(\tau_{A_\al})  \approx \frac{1}{\log (1+Q) - Q/(1+Q)} \brc{|\log \al| -\log (1+Q) + Q +
\sum_{k=1}^\infty \pi_k |\log \pi_k| -1}.
\]

If the prior distribution of the point of change is geometric with a parameter $\rho$,
\[
\pi_k=  \rho(1-\rho)^{k-1}, \quad 0<\rho<1, \quad k \ge 1,
\]
then
\[
\sum_{k=1}^\infty \pi_k |\log \pi_k| = \log \frac{1-\rho}{\rho} - \frac{\log (1-\rho)}{\rho},
\]
and, therefore, the approximation to the average detection delay is given by
\begin{align*}
\ADD^{\pi}(\tau_{A_\al})  &\approx \frac{1}{\log (1+Q) - Q/(1+Q)} \Bl\{|\log \al| -\log (1+Q) + Q
\\
& \quad + \log \frac{1-\rho}{\rho} - \frac{\log (1-\rho)}{\rho} -1\Br\}.
\end{align*}

Note also that in the case of i.i.d.\ observations the detection statistic $G_n$ obeys the
recursion
\[
G_n = (G_{n-1} - \Pi_{n+1}) e^{\De Z_n} + \Pi_{n+1}, \quad G_0=1,
\]
where $\Pi_{n+1}= (1-\rho)^n$ for the geometric prior distribution.

\subsection{Detection of a change in the mean of a Gaussian autoregressive process} \label{ss:AR}

Let $X_n = \te \Ind{\la \le n} + V_n$, $n \ge 1$, where $\te\neq 0$ is a constant ``signal" that
appears at an unknown point in time $\la$ and $V_n$, $n \ge 1$ is zero-mean stable Gaussian $p$-th
order autoregressive process (``noise") $\text{AR}(p)$ that obeys the recursive relation
\[
V_n = \sum_{j=1}^p \de_j V_{n-j} + \xi_n, \quad n \ge 1, \quad \text{$V_j =0$ for $j \le 0$},
\]
where $\xi_n$, $n\ge 1$ are i.i.d. $\Nc(0,\si^2)$ and $1-\sum_{j=1}^p \de_j y^j =0$ has no roots
inside the unit circle.

For $i \ge 1$, define
\[
\tiX_i= \begin{cases}
 X_1 & \text{if $i=1$} \\
 X_i - \sum_{j=1}^{i-1} \de_j X_{i-j} & \text{if $2 \le i \le p$}\\
 X_i - \sum_{j=1}^{p} \de_j X_{i-j} & \text{if $i \ge p+1$},
\end{cases}
\]
and for $i \ge k$ and $k=1,2,\dots$, define
\[
\tite_i= \begin{cases}
 \te & \text{if $i=k$} \\
 \te(1 - \sum_{j=1}^{i-k} \de_j) & \text{if $k+1 \le i \le p+k-1$} \; .\\
 \te(1 - \sum_{j=1}^{p} \de_j) & \text{if $i \ge p+k$}
\end{cases}
\]
The conditional pre-change pdf $f_{0}(X_i|\Xb_1^{i-1})$ is of the form
\begin{align*}
f_{0}(X_i\mid \Xb_1^{i-1}) = \tfrac{1}{\si}\varphi\brc{\tfrac{\tiX_i}{\si}} \quad \text{for all $i
\ge 1$},
\end{align*}
and the conditional post-change pdf $f_{1}(X_i|\Xb_1^{i-1})$, conditioned on $\la=k$, is given by
\begin{align*}
f_{1}(X_i\mid \Xb_1^{i-1})= \tfrac{1}{\si}\varphi\brc{\tfrac{\tiX_i- \tite_{i}}{\si}} \quad
\text{for $i \ge k$},
\end{align*}
where $\varphi(y) = (2\pi)^{-1/2}\exp\set{-y^2/2}$ is the standard normal pdf.

Using these formulas, we easily obtain that the LLR
\begin{equation*}
Z_n^k = \frac{1}{\si^2}\sum_{i=k}^{n} \tite_{i}\tiX_{i}- \frac{1}{2\si^2}\sum_{i=k}^n \tite_{i}^2,
\quad 1\le k \le n, \quad n=1,2,\dots
\end{equation*}
Write
\[
q = \frac{\te^2}{2\si^2} \brc{1-\sum_{j=1}^p\de_j}^2.
\]
Note that, under $\Pb_k$, the LLR process $Z_{n+k-1}^k$, $n \ge 1$ has independent Gaussian
increments $\De Z_n$. Moreover, the increments are i.i.d.\ for $n \ge p+1$ with mean $\Eb_k \De Z_n
= q$ and variance $q/2$. Using this property, it can be shown that $Z_{n+k-1}^k/n$ converges
$r$-quickly to $q$ for all positive $r$ under $\Pb_k$ (see Tartakovsky and Veeravalli
\cite{TartakovskyVeeravalli-TVP2004} for further details and generalizations).

Therefore, Theorem \ref{th4} and Corollary \ref{cor2} can be applied to show that the detection
test $\tau_{A_\al}$ with $A_\al=1/\al$ asymptotically minimizes all positive moments of the
detection delay.

Note also that in the ``stationary" mode when the stopping time $\tau_A \gg k$, the original
problem of detecting a change of the intensity $\te$ in a correlated Gaussian noise is equivalent
to detecting a change of the intensity $\te(1-\sum_{j=1}^p\de_j)$ in white Gaussian noise. This is
primarily because the original problem allows for whitening without loss of information through the
innovations $\tiX_n$, $n \ge 1$ that contain the same information about the hypotheses $H_k$ and
$H_\infty$ as the original sequence $X_n$, $n \ge 1$.

\subsection{Detection of additive changes in state-space hidden Markov models} \label{ss:HMM}

Consider the linear state-space hidden Markov model where the unobserved $m$-dimensional Markov
component $\te_n$ is given by the recursion
\[
\te_n = \Fb \te_{n-1} + W_{n-1} + \nub_\te \Ind{\la \le n}, \quad n \ge 0, \quad \te_0=0,
\]
and the observed $r$-dimensional component
\[
X_n = \te_n + V_n + \nub_x \Ind{\la \le n} , \quad n \ge 1.
\]
Here $W_n$ and $V_n$ are zero-mean Gaussian i.i.d.\ vectors having covariance matrices $\mb{K}_W$
and $\mb{K}_V$, respectively; $\nub_\te=(\nu_\te^1,\dots,\nu_\te^m)$ and
$\nub_x=(\nu_x^1,\dots,\nu_x^r)$ are vectors of the corresponding change intensities; and $\Fb$ is
a $m\times m$ matrix.

It can be shown that under the no-change hypothesis the observed sequence $X_n$, $n \ge 1$ has an
equivalent representation with respect to the innovative process $\xi_n=X_n -
\Eb(\te_n|\Fc_{n-1})$, $n \ge 1$:
\[
X_n = \hat{\te}_n + \xi_n, \quad n \ge 1,
\]
where $\xi_n \sim \Nc(0,\Sib_n)$, $n =1,2,\dots$ are independent Gaussian vectors and
$\hat{\te}_n=\Eb(\te_n|\Fc_{n-1})$ (cf., e.g., Tartakovsky \cite{Tarbook91}). Note that
$\hat{\te}_n$ is the optimal (in the mean-square sense) one-step ahead predictor, i.e., the
estimate of $\te_n$ based on $\Xb_1^{n-1}$, which can be obtained by the Kalman filter. Under the
hypothesis ``$H_k: \la=k$",
\[
X_n = \de_n(k)+ \hat{\te}_n + \xi_n, \quad n \ge 1,
\]
where $\de_n(k)$ depends on $n$ and the change point $k$. The value of $\de_n(k)$ can be computed
using relations given, e.g., in Basseville and Nikiforov \cite{BassevilleNikiforovbook93}.

It follows that the LLR $Z_n^k$ is given by
\[
Z_n^k = \sum_{i=k}^n \de_i(k)^T \Sib_i^{-1} \xi_i - \frac{1}{2} \sum_{i=k}^n \de_i(k)^T \Sib_i^{-1}
\de_i(k),
\]
where $\Sib_i$ are given by Kalman equations (see, e.g., (3.2.20) in
\cite{BassevilleNikiforovbook93}). Therefore, the original abrupt change detection problem that
occurs at $\la=k$ is equivalent to detecting a gradual change from zero to $\de_i(k)$, $i \ge k$ in
the sequence of independent Gaussian innovations $\xi_i$ with the covariance matrices $\Sib_i$.
These innovations can be formed by the Kalman filter. Note also that since the post-change
distribution depends on the change point $k$ through the value of $\de_n(k)$, there is no efficient
recursive formula for the statistic $G_n$ as in the i.i.d.\ case.

As $n \to\infty$, the normalized LLR $n^{-1} Z_{k+n-1}^k$ converges almost surely under $\Pb_k$ to
the positive constant
\[
q= \frac{1}{2} \lim_{n\to\infty} \frac{1}{n}\sum_{i=k}^{k+n-1} \de_i(k)^T \Sib_i^{-1} \de_i(k) .
\]
Using \cite{BassevilleNikiforovbook93}, we obtain that this constant is given by
\[
q=\frac{1}{2}\set{(z\Ibo_m - \Fb^*)^{-1} \nub_\te +[\Ibo_r-(z\Ibo_m-\Fb^*)^{-1}\Fb \mb{K}]\nub_x},
\]
where $\mb{K}$ is the gain in the Kalman filter in the stationary regime, $\Ibo_m$ is the unit
$m\times m$ matrix, and $\Fb^*=\Fb(\Ibo_m-\mb{K})$.

Moreover, since the process $Z_{k+n-1}^k$, $n\ge 1$ is Gaussian with independent increments,
$n^{-1} Z_{k+n-1}^k$ converges strongly completely  to $q$ (i.e., $r$-quickly for all $r>0$, see
Tartakovsky \cite{TarSISP98}). Therefore, Corollary \ref{cor2} shows that the detection test
$\tau_{A_\al}$ is asymptotically optimal as $\alo$ with respect to all positive moments of the
detection delay.

\subsection{Detection of non-additive changes in mixture and HMM models} \label{ss:Mixtures}

In the previous two examples the changes were additive. Consider now an example with non-additive
changes where the observations are i.i.d.\ in the ``out-of-control" mode and mixture-type dependent
in the ``in-control" mode. This example was used by Mei \cite{MeiThesis} as a counterexample to
disprove that the CUSUM and SRP detection tests are asymptotically optimal in the minimax setting
with the lower bound on the mean time to false alarm. However, we show below that the proposed
Bayesian test is asymptotically optimal. This primarily happens because the strong law of large
numbers still holds for the problem considered, while a stronger essential supremum condition (cf.\
Lai \cite{LaiIEEE98}), which is required for obtaining a lower bound for the minimax average
detection delay, fails.

Let $g_1(X_n)$, $g_2(X_n)$, and $f_1(X_n)$ be three distinct densities. The problem is to detect
the change from the mixture density
\[
f_0(\Xb_1^n) = \beta \prod_{i=1}^n g_1(X_i) + (1-\beta) \prod_{i=1}^n g_2(X_i)
\]
to the density $f_1$, where $0 < \beta < 1$ is a mixing probability. Therefore, the observations
are dependent with the joint pdf $f_0(\Xb_1^n)$ before the change occurs and i.i.d.\ with the
density $f_1$ after the change occurs.

Denote $R_j(n)= \log [f_1(X_n)/g_j(X_n)]$ and $I_j = \Eb_{1}R_j(1)$, $j=1,2$.

It is easy to show that
\[
\frac{f_1(X_i)}{f_0(X_i|\Xb_1^{i-1})} = \frac{e^{R_2(i)}(\beta \xi_{i-1} +1 -\beta)}{\beta \xi_{i}
+1 - \beta},
\]
where $\xi_i= \prod_{m=1}^i \De \xi_m$, $\De\xi_m=g_1(X_m)/g_2(X_m)$. Next, note that
\[
\prod_{i=k}^n \frac{1 -\beta + \beta \xi_{i-1}}{1 -\beta + \beta \xi_{i}} = \frac{1+ v
\xi_{k-1}}{1+ v \xi_{n}},
\]
where $v=\beta/(1-\beta)$, so that the LLR
\begin{equation} \label{LLRmix}
Z_n^k := \sum_{i=k}^n \log \frac{f_1(X_i)}{f_0(X_i|\Xb_1^{i-1})}  = \sum_{i=k}^n R_2(i) + \log
\frac{1+ v \xi_{k-1}}{1+ v \xi_{n}}.
\end{equation}

Assume that $I_1 > I_2$, in which case the expectation $\Eb_{k} \log \De\xi_m<0$ for $k <m$ and,
hence,
\[
\xi_n =  \xi_{k-1} \prod_{m=k}^n \De \xi_m \xra[n\to\infty]{\Pb_k-\text{a.s.}} 0 \quad \text{for
every $k<\infty$}.
\]

The condition \eqref{inprobk}, which is necessary for the lower bound \eqref{Dmlowergeninf} to be
satisfied, holds with the constant $q=I_2$. Indeed, since $R_2(i)$, $i \ge k$ are i.i.d.\ random
variables under $\Pb_k$ with mean $I_2$ and since $\xi_n \to 0$, the LLR obeys the strong law of
large numbers:
\[
\frac{1}{n} Z^k_{n+k-1} \to I_2 \quad \text{$\Pb_k$-a.s. as $n \to\infty$},
\]
which implies \eqref{inprobk} with $q=I_2$ and, hence, the lower bound \eqref{Dmlowergeninf},
\[
\inf_{\tau\in\Dal} \D_m^\pi(\tau) \ge \brc{\frac{|\log\al|}{I_2}}^m (1+o(1)) \quad \text{as $\alo$
for all $m>0$}.
\]

Next, using \eqref{logGnnew} and \eqref{LLRmix}, we can write the statistic $\log G_n$ is the
following form
\[
\log G_n = \sum_{i=k}^n R_2(i) + \psi(k,n) + \log [\pi_k(1+Y_n^k)],
\]
where
\[
\psi(k,n)=\log \frac{1+ v \xi_{k-1}}{1+ v \xi_{n}}.
\]
The sequence $Y_n^k$, $n \ge k$ is slowly changing by the argument given in the proof of Theorem
\ref{th3}. The sequence $\psi(k,n)$, $n \ge k$ is also slowly changing. In fact, since $\xi_n \to
0$ w.p.\ 1, it converges to the finite random variable $\log(1+v\xi_{k-1})$. Therefore, by the
nonlinear renewal theorem \cite{Woodroofebook82},
\[
\ADD(\tau_{A_\al}) = \brc{\frac{|\log \al| }{I_2}}(1+o(1)) \quad \text{as $\alo$},
\]
and the detection procedure $\tau_{A_\al}$ is asymptotically optimal.

Note that the results of Tartakovsky and Veeravalli \cite{TartakovskyVeeravalli-TVP2004} suggest
that the Shiryaev detection procedure is also asymptotically optimal under the traditional
constraint on the average false alarm probability. On the other hand, as we mentioned above, the
minimax property of the CUSUM and Shiryaev-Roberts tests does not hold in the example considered.

Finally, we note the above simple mixture model is obviously a degenerate case of a more general
model governed by a two-state HMM when transition probabilities between states are equal to zero
and the initial distribution is given by the probability $\beta$. The proposed Bayesian procedure
(as well as the Shiryaev procedure in the conventional setting) remains asymptotically optimal for
the model where the pre-change distribution is controlled by a finite-state (non-degenerate) HMM,
while the post-change model is i.i.d. On the other hand, the condition C1 of Fuh \cite{Fuh03} does
not hold and, therefore, one may not conclude that the CUSUM test is minimax asymptotically optimal
under the constraint on the average run length to false alarm. For such a model, the minimax
asymptotic optimality property of the CUSUM is an open problem. Simulation results show that the
performance of the CUSUM test is poor at least for the moderate false alarm rate, while the
performance of the Bayesian tests is high. Further details will be presented elsewhere.

\section{Concluding Remarks} \label{s:Remarks}

1. As we already mentioned in the introduction, the global false alarm probability constraint
$\sup_k\Pb_k(\tau<k) = \Pb_\infty(\tau<\infty)\le\al$ leads to an unbounded worst-case expected
detection delay $\sup_k \Eb_k(\tau-k|\tau \ge k)$ whenever $\al<1$ due to a high price that should
be paid for such a strong constraint.  Note that to overcome this difficulty in a minimax setting a
dynamic sampling technique can be used when it is feasible (cf.\ Assaf et al \cite{AssafetalAS93}).
To the expense of a large amount of data that must be sampled, the worst-case average detection
delay may then be made bounded, yet keeping the global PFA below the given small level. However,
dynamic sampling is rarely possible in applications. We, therefore, considered a Bayesian problem
with the prior distribution. The proposed asymptotically Bayesian detection test can be regarded as
the Shiryaev detection procedure with a threshold that increases in time. The need for the
threshold increase is due to the strong constraint imposed on the global PFA in place of the
average PFA constraint used in Shiryaev's classical problem setting.

2. While the results of the present paper may be used to devise a reasonably simple detection
procedure to handle the global probability bound on false alarms, the author's personal opinion is
that this constraint is too strong to be useful in applications. In fact, the conditional ADD
$\Eb_k(\tau_A-k|\tau_A\ge k)$ of the proposed detection procedure grows fairly fast with $k$, and
the ``nice" property that the Bayesian ADD is as small as possible (for small $\al$) perhaps will
not convince practitioners in the usefulness of the test. In addition, the mean time to false alarm
in this detection procedure is unbounded, which is an unavoidable recompense for the very strong
global PFA constraint.

3. Taking into account the previous remark, we argue that imposing the bound on the local PFA
$\sup_k \Pb_\infty(k \le \tau \le k +T-1)$ or on the local conditional PFA $\sup_k \Pb_\infty(k \le
\tau \le k +T-1| \tau \ge k)$ is a much more practical approach. The latter conditional PFA is
indeed a proper measure of false alarms in a variety of surveillance problems, as was discussed in
Tartakovsky \cite{TartakovskyIEEECDC05}. It can be then shown that the conventional CUSUM and SRP
detection tests are optimal in the minimax sense for any time window $T$, and asymptotically
uniformly optimal (i.e., for all $k\ge 1$) if the size of the window $T$ goes to infinity at a
certain rate (cf.\ Lai \cite{LaiJRSS95,LaiIEEE98} and Tartakovsky \cite{TartakovskyIEEECDC05}).

4. The sufficient conditions for asymptotic optimality postulated in Theorems \ref{th1} and
\ref{th4} are quite general and hold in most applications. We verified these conditions for the
three examples that cover both additive and non-additive changes in non-i.i.d.\ models. While we
are not aware of the non-i.i.d.\ models reasonable for practical applications for which these
conditions do not hold, such examples may still exist. However, we believe that such situations
should be handled on a case by case basis.

5. Similar results can be proved for general continuous-time stochastic models. A proof of the
lower bound for moments of the detection delay is absolutely identical to the proof of
Theorem~\ref{th1}. However, derivation of the upper bound is not straightforward and requires
certain additional conditions analogous to those used in Baron and Tartakovsky
\cite{BaronTartakovskySA05}.

\section*{Acknowledgements}

This research was supported in part by the U.S.\ Office of Naval Research grant N00014-06-1-0110 at
the University of Southern California and by the U.S.\ ARMY SBIR contract W911QX-04-C-0001 at
ADSANTEC.

\section*{Appendix}

\renewcommand{\proofname}{\bf Proof of Lemma~\ref{le2}.}

\begin{proof} 

Define $\ga_{\vae,\al}^{(k)}(\tau) =\Pb_k\set{k \le \tau < k+(1-\vae)L_\al}$, where
$L_\al=q^{-1}|\log \al|$. A quite tedious argument analogous to that used in the proof of Lemma 1
of Tartakovsky and Veeravalli ~\cite{TartakovskyVeeravalli-TVP2004} yields
\begin{align*} 
\ga_{\vae,\al}^{(k)}(\tau) & \le e^{-(1-\vae^2)\log \al}
\Pb_{\infty}\set{k\le\tau<k+ (1-\vae)L_\al}
\\
& \quad + \Pb_k\set{\max_{0\le n<(1-\vae)L_\al} Z_{k+n}^k\ge(1-\vae^2)q L_\al}.
\end{align*}
Since $\Pb_{\infty}\set{k\le\tau<k+ (1-\vae)L_\al} \le \Pb_\infty (\tau < \infty)
\le \al$ for any $\tau \in \Deb_\infty(\al)$, we obtain
\begin{equation} \label{gakupper}
\ga_{\vae,\al}^{(k)}(\tau) \le \al^{\vae^2} + \beta_k(\al,\vae),
\end{equation}
where
\[
\beta_k(\al,\vae)=\Pb_k\set{\max_{1\le n \le (1-\vae)L_\al} Z_{k+n-1}^k\ge(1-\vae^2)q L_\al}.
\]

Let $N_\al=\lfloor\vae L_\al\rfloor$ be the greatest integer number $\le \vae
L_\al$. Evidently,
\begin{equation*} 
\ga_{\vae,\al}^\pi(\tau)=\sum_{k=1}^\infty\pi_k\ga_{\vae,\al}^{(k)}(\tau)\le
\sum_{k=1}^{N_\al}\pi_k\ga_{\vae,\al}^{(k)}(\tau) + \Pi_{N_\al+1}
\end{equation*}
and, therefore,
\begin{equation} \label{gammaineq}
\ga_{\vae,\al}^\pi(\tau) \le \Pi_{N_\al+1} + \al^{\vae^2} + \sum_{k=1}^{N_\al}\pi_k
\beta_k(\al,\vae).
\end{equation}
The first two terms go to 0 as $\al\to 0$ for any $\vae>0$. The third term goes to
zero as $\al\to 0$ by condition \eqref{inprobk} and Lebesgue's dominated convergence
theorem. Since the right side in \eqref{gammaineq} does not depend on $\tau$, this
completes the proof of \eqref{ga0general}.

Using the inequality $\Pb_\infty (\tau_A < \infty) \le 1/A$ and applying the same
argument as above shows that
\begin{equation} \label{gammaineqtauA}
\ga_{\vae,A}^\pi(\tau_A) \le \Pi_{N_A+1} + 1/A^{\vae^2} + \sum_{k=1}^{N_A}\pi_k \beta_k(A,\vae),
\end{equation}
where $N_A=\lfloor\vae L_A\rfloor$ and
\[
\beta_k(A,\vae)=\Pb_k\set{\max_{1\le n \le (1-\vae)q^{-1} \log A} Z_{k+n-1}^k\ge(1-\vae^2) \log A}.
\]
Again all three terms on the right-hand side of \eqref{gammaineqtauA} tend to zero
as $A\to \infty$, which proves \eqref{ga0tauA}.

\end{proof}


\end{document}